\documentclass[reqno,12pt,a4]{amsart}
\usepackage{amssymb}
\usepackage{amsmath}

\usepackage[all]{xy}

\def\Hom{\mathop{\rm Hom}\nolimits}
\def\Ext{\mathop{\rm Ext}\nolimits}

\def\Ad{\mathop{\rm Ad}\nolimits}

\def\ev{\mathop{\rm ev}\nolimits}

\newtheorem{thm}{Theorem}[section]
\newtheorem{lem}[thm]{Lemma}

\newtheorem{cor}[thm]{Corollary}

\newtheorem{rmk}[thm]{Remark}

\numberwithin{equation}{section}

\newcommand{\cone}{\operatorname{cone}}

\date{\today}
\author{V. Manuilov and K. Thomsen}
\title[Translation invariant asymptotic homomorphisms]{Extensions of $C^*$-algebras and translation invariant asymptotic homomorphisms}

\begin{document}
\maketitle

\begin{abstract}
Let $A$, $B$ be $C^*$-algebras; $A$ separable, $B$ $\sigma$-unital and stable.
We introduce a notion of translation invariance for asymptotic homomorphisms from $SA=C_0(\mathbb R)\otimes A$ to $B$ and show that the Connes--Higson construction applied to any extension of $A$ by $B$ is homotopic to a translation invariant asymptotic homomorphism. In the other direction we give a construction which produces extensions of $A$ by $B$ out of such a translation invariant asymptotic homomorphism. This leads to our main result; that the homotopy classes of extensions coincide with the homotopy classes of translation invariant asymptotic homomorphisms.

\end{abstract}

\section{Introduction}

The excision properties of the E-theory of Connes and Higson hinges on a fundamental construction, introduced in \cite{CH}, which associates an asymptotic homomorphism $SA \to B$ to an extension of the separable $C^*$-algebra $A$ by another $\sigma$-unital $C^*$-algebra $B$. This construction transforms an extension into something which seems to be of a different nature, in much the same way as the Busby-invariant transforms extensions, viewed as short exact sequences of $C^*$-algebras, into $*$-homomorphisms. Both constructions have turned out to be very useful for the study of extensions of $C^*$-algebras, and in this way also for many other purposes. But while it is easy to see that the Busby-invariant holds all relevant information on the extension, it is far more difficult to decide how much information is lost when considering the asymptotic homomorphism arising from an extension, rather than the extension itself. This question is quite intriguing because the approach to the study of extensions which is based solely on the Busby-invariant, i.e. the approach developed by Brown, Douglas and Fillmore in \cite{BDF2}, and later in higher generality by Kasparov \cite{K}, is now known to have serious deficiencies: It may not have a group structure, \cite{An}, not even in cases where the homotopy classes of extensions \emph{do} have such a structure, \cite{Kirch}. In our previous work, \cite{MT3,MT4}, we have been able, in certain cases, to determine exactly which structure the asymptotic homomorphism arising from the Connes-Higson construction encodes. But the general case has escaped us, and the results of the present paper show why. As we shall explain the asymptotic homomorphism $SA \to B$ has a property which has been overlooked so far; it is namely \emph{translation invariant} in the sense that a shift in the parameter of the asymptotic homomorphism corresponds exactly to a shift of the same amount in $SA$. This observation shows that the parameter of the asymptotic homomorphism is not just a structureless dummy, whose only raison d'\^etre is to approach infinity without actually getting there, but that it has a serious relationship to the $C^*$-algebra $SA$ on which the asymptotic homomorphism is defined, at least when it arises from the Connes-Higson construction. Besides making this observation precise, we will show that this additional structure of the asymptotic homomorphism is exactly what is needed in order to recover the extension, showing that there are no more surprises to be uncovered.

Let $A$ and $B$ be $C^*$-algebras. We say that a one-parameter family of maps, $\psi = (\psi_t)_{t \in \mathbb R} : A \to B$, is an \emph{asymptotic homomorphism} when $(\psi_t)_{t \in [1,\infty)}$ is (in the usual sense, \cite{CH}) and $\lim_{t\to-\infty}\psi(a)=0$ for all $a\in A$. $\psi$ is \emph{uniformly continuous} when the function $t \mapsto \psi_t(a)$ is uniformly continuous for all $a \in A$, and \emph{equi-continuous} when the family of maps $\psi_t : A \to B, t \in \mathbb R$, is. 

Let $\tau=(\tau_s)_{s\in\mathbb R}$ be the action of $\mathbb R$ on $SA = C_0(\mathbb R, A)$ by translations,
$$
\tau_s(f)(t) = f(t -s),\quad f\in SA,
$$ 
and consider an asymptotic homomorphism $\psi = (\psi_t)_{t \in \mathbb R} : SA \to B$. We say that $\psi$ is \emph{asymptotically translation invariant} when $\lim_{t \to \infty} \psi_{t-s}(f) - \psi_t \circ \tau_s(f) =0$ for all $s \in \mathbb R$ and all $f \in SA$, and \emph{translation invariant} when $\psi$ is equicontinuous, and
$$
\psi_t \circ \tau_s = \psi_{t-s}
$$
for all $t,s \in \mathbb R$. Two (asymptotically) translation invariant asymptotic homomorphisms, $\varphi$ and $\psi$, are \emph{homotopic} when there is a (asymptotically) translation invariant asymptotic homomorphism $\Phi : SA \to IB (=C[0,1]\otimes B)$ such that $\ev_1 \circ \Phi_t = \varphi_t$ and $\ev_0 \circ \Phi_t = \psi_t$ for all $t \in \mathbb R$, where $\ev_s : IB \to B$ denotes evaluation at $s \in [0,1]$. We denote by $[[SA,B]]_{\tau}$ the set of homotopy classes of translation invariant asymptotic homomorphisms $SA \to B$, and by $[[SA,B]]_{a,\tau}$ the set of homotopy classes of asymptotically translation invariant asymptotic homomorphisms $SA \to B$. There is then a natural map
$$
[[SA,B]]_{\tau} \to [[SA,B]]_{a,\tau}.
$$
It will be a corollary of our main result that this map is an isomorphism of abelian semigroups when $A$ is separable and $B$ $\sigma$-unital and stable.

\section{From extensions to translation invariant asymptotic homomorphisms}\label{vch}

Let $A$ and $B$ be $C^*$-algebras; $A$ separable, $B$ $\sigma$-unital. The multiplier algebra of $B$ will be denoted by $M(B)$ and the generalized Calkin algebra $M(B)/B$ by $Q(B)$. Let $\psi : A \to Q(B)$ be an extension. We choose
\begin{enumerate}
\item[a)] a strictly positive element $b_0 \in B$ such that $0 \leq b_0 \leq 1$,
\item[b)] a continuous and homogeneous lift $s : A \to M(B)$ of $\psi$, and
\item[c)] a sequence $h_0 \leq h_1 \leq h_2 \leq \dots$ of functions in $C_0(0,1]$,
\end{enumerate}
such that 
\begin{enumerate}
\item[d)] $h_{n+1}h_n = h_n$ for all $n$,
\end{enumerate}
and the sequence $u_n = h_n(b_0), n = 0,1,2, \dots$, has the properties
\begin{enumerate}
\item[i)] $ \lim_{n \to \infty} u_ns(a) - s(a)u_n = 0$ for all $a \in A$, and
\item[ii)] $\lim_{n \to \infty} u_nb =b$ for all $b \in B$.
\end{enumerate}
This is possible by the arguments that prove the existence of quasi-central approximate units, \cite{A}. $\{u_n\}_{n=0}^{\infty}$ is a \emph{unit sequence} in the sense of \cite{MT3} and \cite{MT5}.

We set $\Delta_0 = \sqrt{u_0}$ and $\Delta_n = \sqrt{u_n - u_{n-1}}, n \geq 1$. A crucial observation is that thanks to i) we have that
\begin{equation}\label{e13}
\Delta_i\Delta_j =0, \ \text{when $|i-j| \geq 2$.}
\end{equation}
Let $t_0 \leq t_1 \leq t_2 \leq \dots$ be a sequence in $[0, \infty[$ such that
\begin{enumerate}
\item[iii)] $\lim_{n \to \infty} t_n = \infty$,
\item[iv)] $\lim_{n \to \infty} t_n - t_{n+1} = 0$.
\end{enumerate}

The following lemma appears in \cite{MT5}.

\begin{lem}\label{mt5} For any norm-bounded sequence $\{m_j\} \subseteq M(B)$, and for any $k\in\mathbb N$, the sum $\sum_{j=0}^{\infty} \Delta_j m_j
\Delta_{j+k}$ converges in the strict topology to an element of $M(B)$, and
$$
\left\|\sum_{j=0}^{\infty} \Delta_j m_j \Delta_{j+k} \right\| \leq \sup_j \|m_j\|.
$$
Furthermore, $\sum_{j=0}^{\infty} \Delta_j m_j \Delta_{j+k} \in B$ when $\lim_{j \to \infty} \|m_j\| = 0$.
\end{lem}

For each $t \in \mathbb R$, define $\varphi_t : SA \to B$ by
\begin{equation}\label{ugh}
\varphi_t(f)= \sum_{j=0}^{\infty} \Delta_j s\left(f(t -t_j)\right)\Delta_j .
\end{equation}
The sum converges in norm because $\lim_{j \to \infty} s\left(f(t -t_j)\right) = 0$, cf. Lemma \ref{mt5}.

Since $s$ is homogeneous and continuous at $0$, there is an $L > 0$ such that $\|s(a)\| \leq L\|a\|$ for all $a \in A$. It follows that
\begin{equation}
\left\|\varphi_t(f)\right\| \leq L\|f\|
\end{equation}
for all $f$ and $t$, cf. Lemma \ref{mt5}.

\begin{lem}\label{a1}

$\left(\varphi_t\right)_{t \in \mathbb R}$ is a translation invariant asymptotic homomorphism. I.e.

\begin{enumerate}
\item[A)] $\left(\varphi_t\right)_{t \in [0,\infty)}$ is an asymptotic homomorphism,
\item[B)] $\varphi_t \circ \tau_s = \varphi_{t-s}$ for all $t,s \in \mathbb R$,
\item[C)] $\lim_{t \to -\infty}\varphi_t(f) = 0$ for all $f \in SA$,
\item[D)] $\varphi$ is equicontinuous.
\end{enumerate}

\end{lem}
\begin{proof} B) is a trivial observation and C) follows from Lemma \ref{mt5} since $\lim_{x \to -\infty} s(f(x)) = 0$, so it remains to prove A) and D). We begin with A): Let $f,g\in SA$ and let $\epsilon > 0$. By uniform continuity of the map $x \mapsto s(g(x))$, there is a $\delta > 0$ so small that $\left\|s(g(x)) - s(g(y))\right\| < \epsilon$ when $|x-y| < \delta$. Choose $N \in \mathbb N$ so large that $t_j - t_{j-1} < \delta$ when $j \geq N$. Since $\lim_{t \to \infty} s(k(t-t_j)) = 0$ for all $j \in \mathbb N$ and all $k \in SA$, there is then a $T \in \mathbb N$ so large that
$$
\left\|\varphi_t(f) -  \sum_{j=N}^{\infty} \Delta_j s\left(f(t -t_j)\right)\Delta_j \right\| < \epsilon,
\qquad
\left\|\varphi_t(g) -  \sum_{j=N}^{\infty} \Delta_j s\left(g(t -t_j)\right)\Delta_j \right\| < \epsilon,
$$
and
$$
\left\|\varphi_t(fg) -  \sum_{j=N}^{\infty} \Delta_j s\left(fg(t -t_j)\right)\Delta_j \right\| < \epsilon
$$
when $t \geq T$.
In what follows we use the notation $a\sim_\epsilon b$ when $\|a-b\|\leq\epsilon$.
Then
\begin{align*}
\varphi_t(f)\varphi_t(g)
&\sim_{(\|f\| + \|g\|)L\epsilon} \\ \sum_{j=N}^{\infty} \Delta_j s\left(f(t-t_j)\right)\Delta_j^2s\left(g(t - t_j)\right)&\Delta_j +   \sum_{j=N}^{\infty} \Delta_j s\left(f(t-t_j)\right)\Delta_j\Delta_{j+1}s\left(g(t - t_{j+1})\right)\Delta_{j+1}   \\ & + \sum_{j=N+1}^{\infty} \Delta_j s\left(f(t-t_j)\right)\Delta_j\Delta_{j-1}s\left(g(t - t_{j-1})\right)\Delta_{j-1}\\
\end{align*}
when $t \geq T$. It follows from i) that if $N$ is large enough, then
$$
\sup_{x \in \mathbb R} \left\|\left[\Delta_j^2,s\left(g(x))\right)\right]\right\| < \epsilon ,
\qquad
\sup_{x \in \mathbb R} \left\|\left[\Delta_j\Delta_{j\pm 1},s\left(g(x)\right)\right]\right\| < \epsilon,
$$
when $j \geq N$. So by increasing $N$ if necessary, we may assume, thanks to Lemma \ref{mt5}, that
\begin{align*}
 \sum_{j=N}^{\infty} \Delta_j s\left(f(t-t_j)\right)\Delta_j^2s\left(g(t - t_j)\right)\Delta_j +   \sum_{j=N}^{\infty} \Delta_j s\left(f(t-t_j)\right)\Delta_j\Delta_{j+1}s\left(g(t - t_{j+1})\right)\Delta_{j+1}   \\ + \sum_{j=N+1}^{\infty} \Delta_j s\left(f(t-t_j)\right)\Delta_j\Delta_{j-1}s\left(g(t - t_{j-1})\right)\Delta_{j-1}\\  \\
\sim_{\epsilon} \sum_{j=N}^{\infty} \Delta_j s\left(f(t-t_j)\right)s\left(g(t - t_j)\right)\Delta_j^3 +   \sum_{j=N}^{\infty} \Delta_j s\left(f(t-t_j)\right)s\left(g(t - t_{j+1})\right)\Delta_j\Delta_{j+1}^2   \\ + \sum_{j=N+1}^{\infty} \Delta_j s\left(f(t-t_j)\right)s\left(g(t - t_{j-1})\right)\Delta_j\Delta_{j-1}^2\\  \\
\end{align*}
when $t \geq T$. It follows from ii) that
$$
\lim_{j \to \infty} \sup_{x \in \mathbb R} \left\|\Delta_j^2 \left[ s(f(x))s(g(x)) - s(f(x)g(x))\right]\right\| = 0
$$
so by increasing $N$ further, and using Lemma \ref{mt5} again, we can arrange that
\begin{align*}
 \sum_{j=N}^{\infty} \Delta_j s\left(f(t-t_j)\right)s\left(g(t - t_j)\right)\Delta_j^3 +   \sum_{j=N}^{\infty} \Delta_j s\left(f(t-t_j)\right)s\left(g(t - t_{j+1})\right)\Delta_j\Delta_{j+1}^2   \\ + \sum_{j=N+1}^{\infty} \Delta_j s\left(f(t-t_j)\right)s\left(g(t - t_{j-1})\right)\Delta_j\Delta_{j-1}^2   \\ \sim_{\epsilon}
\sum_{j=N}^{\infty} \Delta_j s\left(f(t-t_j)g(t - t_j)\right)\Delta_j^3 +   \sum_{j=N}^{\infty} \Delta_j s\left(f(t-t_j)g(t - t_{j+1})\right)\Delta_j\Delta_{j+1}^2   \\ + \sum_{j=N+1}^{\infty} \Delta_j s\left(f(t-t_j)g(t - t_{j-1})\right)\Delta_j\Delta_{j-1}^2  .
\end{align*}
when $t \geq T$. By using iv) we can also arrange that
\begin{align*}
\sum_{j=N}^{\infty} \Delta_j s\left(f(t-t_j)g(t - t_j)\right)\Delta_j^3 +   \sum_{j=N}^{\infty} \Delta_j s\left(f(t-t_j)g(t - t_{j+1})\right)\Delta_j\Delta_{j+1}^2   
\\
 + \sum_{j=N+1}^{\infty} \Delta_j s\left(f(t-t_j)g(t - t_{j-1})\right)\Delta_j\Delta_{j-1}^2  
\\
\sim_{\epsilon}\sum_{j=N}^{\infty} \Delta_j s\left(f(t-t_j)g(t - t_j)\right)\Delta_j^3 +   \sum_{j=N}^{\infty} \Delta_j s\left(f(t-t_j)g(t - t_{j})\right)\Delta_j\Delta_{j+1}^2   \\ + \sum_{j=N+1}^{\infty} \Delta_j s\left(f(t-t_j)g(t - t_{j})\right)\Delta_j\Delta_{j-1}^2 
\end{align*} 
when $t \geq T$. Since $\Delta_j^3 + \Delta_j\Delta_{j+1}^2 + \Delta_j\Delta_{j-1}^2 = \Delta_j$, we can finally arrange that
\begin{align*}
\sum_{j=N}^{\infty} \Delta_j s\left(f(t-t_j)g(t - t_j)\right)\Delta_j^3 +   \sum_{j=N}^{\infty} \Delta_j s\left(f(t-t_j)g(t - t_{j})\right)\Delta_j\Delta_{j+1}^2  
\\
 + \sum_{j=N+1}^{\infty} \Delta_j s\left(f(t-t_j)g(t - t_{j})\right) 
\Delta_j\Delta_{j-1}^2  
\\
\sim_{\epsilon}  \sum_{j=N}^{\infty} \Delta_j s(fg(t -t_j))\Delta_j 
\end{align*} 
when $t \geq T$. All in all it follows that there is a $T \in \mathbb R$ so that
$$
 \varphi_t(f)\varphi_t(g) \sim_{(\|f\| + \|g\|)L\epsilon + 5\epsilon} \varphi_t(fg)
$$
for all $t \geq T$. The asymptotic linearity and self-adjointness of $\varphi$ follows in the same way.

D): Let $\epsilon > 0$. For each $x \in A, \delta > 0$, let $B(x,\delta) = \{a \in A: \ \|a-x\| < \delta\}$. Since $K = f(\mathbb R)$ is a compact subset of $A$, and $s$ is continuous, there are finite sets, $\{x_1,x_2, \dots, x_N\} \subseteq K$ and $\{\delta_1,\delta_2, \dots, \delta_N\} \subseteq ]0,1[$, such that $K \subseteq \bigcup_{i=1}^N B(x_i, \delta_i)$ and $a,b \in B(x_i,2 \delta_i) \Rightarrow \|s(a) - s(b)\| < \epsilon $ for all $i$. Set $\delta = \min \{\delta_1,\delta_2, \dots, \delta_N\}$. It follows then from Lemma \ref{mt5} that when $g \in SA$ and $\|f-g\| < \delta$, then $\left\|\varphi_t(f)-\varphi_t(g)\right\| \leq \epsilon$ for all $t \in \mathbb R$.

\end{proof}

The next aim is to show that up to homotopy of translation invariant asymptotic homomorphisms, $\varphi$ depends only on $\psi$.

\subsection{Independence of the choice of $\{t_n\}_{n=0}^{\infty}$.}

When $\{t_n'\}_{n=0}^{\infty}$ is another sequence in $[0,\infty[$ such that iii) and iv) hold, we define $\Phi_t^{\lambda} : SA \to B$ by
$$
\Phi_t^{\lambda}(f) = \sum_{j=0}^{\infty} \Delta_j s(f(t-\lambda t_j - (1-\lambda)t_j'))\Delta_j .
$$

To see that $[0,1] \ni \lambda \mapsto \Phi_t^{\lambda}(f)$ is continuous, let $\epsilon > 0$ be given. Since $\lim_{j \to \infty} t_j = \lim_{j \to \infty} t_j' = \infty$, there is an $N$ such that $\left\|\Phi_t^{\lambda}(f) - \sum_{j=0}^{N} \Delta_j s(f(t-\lambda t_j - (1-\lambda)t_j'))\Delta_j\right\| < \epsilon$ for all $\lambda \in [0,1]$. Since  $[0,1] \ni \lambda \mapsto  \sum_{j=0}^{N} \Delta_j s(f(t-\lambda t_j - (1-\lambda)t_j'))\Delta_j$ is clearly continuous, we see that $[0,1] \ni \lambda \mapsto \Phi_t^{\lambda}(f)$ is continuous.

We can then define $\Phi_t : SA \to IB (= C[0,1]\otimes B)$ such that
$$
\Phi_t(f)(\lambda) = \Phi_t^{\lambda}(f) =  \sum_{j=0}^{\infty} \Delta_j s(f(t-\lambda t_j - (1-\lambda)t_j'))\Delta_j .
$$
The arguments that proved Lemma \ref{a1} show that $\Phi = \left(\Phi_t\right)_{t \in \mathbb R}$ is a translation invariant asymptotic homomorphism. It clearly defines a homotopy of translation invariant asymptotic homomorphisms connecting the translation invariant asymptotic homomorphisms arising from the two choices of sequences $\{t_n\}_{n=0}^{\infty}$ and $\{t_n'\}_{n=0}^{\infty}$.

\subsection{Independence of the choice of unit sequence.}

Let $g_0 \leq g_1 \leq g_2 \leq \dots$ be another sequence of functions in $C_0(0,1]$ such that d) holds, and $v_n = g_n(b_0), n = 0,1,\dots $, satisfies i)-ii). We can then find a third sequence $k_0 \leq k_1 \leq k_2 \leq \dots$ in $C_0(0,1]$ such that d) holds and the sequence $w_n = k_n(b_0), n = 0,1,2, \dots$, also satisfies i)-ii), and at the same time that
$$
w_nu_n = u_n, \ w_nv_n = v_n
$$
for all $n$. Set $d_n = 1 - \sum_{j=n+1}^{\infty} 2^{-j}, n = 0,1,2, \dots$, and define $W_n \in C[0,1]\otimes B$ such that
$$
W_n(\lambda) = \begin{cases} u_n, & \lambda \in [0,d_n],  \\ \frac{d_{n+1} - \lambda}{d_{n+1} -d_n}u_n + \frac{\lambda - d_n}{d_{n+1} -d_n}w_n, & \lambda \in [d_n,d_{n+1}], \\ w_n, & \lambda \in [d_{n+1},1] . \end{cases}
$$
Define $\Phi_t : SA \to IB$ such that
$$
\Phi_t(f) = \sqrt{W_0}s(f(t-t_0))\sqrt{W_0} + \sum_{j=1}^{\infty} \sqrt{W_j - W_{j-1}}s(f(t-t_j)) \sqrt{W_j - W_{j-1}}.
$$
The arguments of Lemma \ref{a1} show that $\Phi = \left(\Phi_t\right)_{t \in \mathbb R}$ is a translation invariant asymptotic homomorphism, so it gives us a homotopy of translation invariant asymptotic homomorphisms connecting the translation invariant asymptotic homomorphisms arising from the two choices of sequences $\{u_n\}_{n=0}^{\infty}$ and $\{w_n\}_{n=0}^{\infty}$. Using $\{v_n\}_{n=0}^{\infty}$ instead of $\{u_n\}_{n =0}^{\infty}$ in the construction of $\Phi$, we get a homotopy of translation invariant asymptotic homomorphisms connecting the translation invariant asymptotic homomorphisms arising from the two choices of sequences $\{v_n\}_{n=0}^{\infty}$ and $\{w_n\}_{n=0}^{\infty}$. By transitivity of homotopy, we see that the two translation invariant asymptotic homomorphisms arising from $\{v_n\}_{n=0}^{\infty}$ and $\{u_n\}_{n=0}^{\infty}$ are homotopic.

\subsection{Independence of the choice of strictly positive element $b_0$.}\label{subsec}

Let $b_0'$ be another strictly positive element of $B$. Then
$$
c_0(\lambda) = \lambda b_0 + (1 -\lambda)b_0'
$$
is a strictly positive element of $IB$. Constructing an asymptotic homomorphism $SA \to IB$ from $c_0$ as in Lemma \ref{a1}, we obtain a homotopy of translation invariant asymptotic homomorphisms connecting two translation invariant asymptotic homomorphisms $SA \to B$; one arising from the construction by use of $b_0$, the other by use of $b_0'$.

\subsection{Independence of the choice of section $s$.}

When $s' : A \to M(B)$ is another continuous and homogeneous lift of $\psi$ we can choose the sequence $h_0 \leq h_1 \leq h_2 \leq \dots$ such that
$$
\lim_{j \to \infty} \sup_{x \in \mathbb R}\left\|\Delta_j\left(s(f(x)) - s'(f(x))\right)\right\| = 0
$$
for all $f \in SA$. It follows then (from Lemma \ref{mt5}) that
$$
\lim_{t \to \infty} \left[\sum_{j=0}^{\infty} \Delta_j s(f(t-t_j))\Delta_j - \sum_{j=0}^{\infty} \Delta_j s'(f(t-t_j))\Delta_j \right] = 0 
$$
for all $f \in SA$.

\bigskip

\subsection{Homotopy invariance.}

Let $\psi':A\to Q(B)$ be an extension weakly homotopic to $\psi$. Then there exist an extension $\Psi:A\to Q(IB)$ with evaluation maps at the endpoints coinciding with $\psi$ and $\psi'$, respectively. Once more, constructing an asymptotic homomorphism $SA \to IB$ from $\Psi$ as in Lemma \ref{a1} gives us a homotopy of translation invariant asymptotic homomorphisms.

We can now summarize in the following:

\begin{thm}\label{a4} The homotopy class of the translation invariant asymptotic homomorphism $\varphi$ of Lemma \ref{a1} depends only on the weak homotopy class of $\psi$ in the set of extensions of $A$ by $B$.
\end{thm}

In other words, if we denote by $\Ext_h(A,B)$ the (weak) homotopy classes of extensions of $A$ by $B$, the formula (\ref{ugh}) gives rise to a map
$$
CH_{\tau} : \Ext_h(A,B) \to [[SA,B]]_{\tau}.
$$

\section{Relation to the Connes-Higson construction.}

In this section we show that up to homotopy $CH_{\tau}(\psi)$ agrees with the Connes-Higson construction $CH(\psi)$, as it was introduced in \cite{CH}.

Let
$$
\cone = \{f \in C_b(\mathbb R): \ \lim_{t \to - \infty} f(t) = 0, \ \text{and} \ \lim_{t \to \infty} f(t) \ \text{exists} \ \},
$$
and
$$
I = \{f \in C_b(\mathbb R): \ \lim_{t \to \pm \infty} f(t)  \ \text{exist} \ \ \}.
$$
Then the sum
$$
\sum_{j=0}^{\infty} \Delta_j^2 h(t-t_j)
$$
converges in the strict topology to an element $\lambda_t(h) \in M(B)$ for all $t \in \mathbb R$ and every $h \in I$, cf. Lemma \ref{mt5}.

Note that $\lambda_t : I \to M(B)$ is a linear completely positive contraction for all $t$. We claim that
\begin{equation}\label{20}
\lim_{t \to \infty} \lambda_t(f)\lambda_t(g) - \lambda_t(fg) = 0
\end{equation}
for all $f,g \in I$. To see this, observe that
\begin{eqnarray*}
&& \Bigl(\sum_{j=0}^{\infty} \Delta_j^2 f(t-t_j) \Bigr) \Bigl(\sum_{j=0}^{\infty} \Delta_j^2g(t-t_j)\Bigr) \\
&& = \sum_{j=0}^{\infty} \Delta_j^4 f(t-t_j )g(t-t_j) +  \sum_{j=0}^{\infty} \Delta_j^2\Delta_{j+1}^2f(t-t_j )g(t-t_{j+1}) \\
&& \ \ \ \ \ \ \ \ \ \ \ \ \ \ \ \ \ \ \ \ \ \ \ \ \ +
\sum_{j=1}^{\infty} \Delta_{j}^2\Delta_{j-1}^2 f(t-t_{j} )g(t-t_{j-1})
\ \ \ \ \ \ \ \ \ \ \ \ \ \ \text{(by (\ref{17}))} \\
\end{eqnarray*}
If $t$ is big enough,
$$
\sup_{j \geq 0} |f(t-t_{j} )g(t-t_{j+1}) -  f(t-t_{j} )g(t-t_j)| \leq \epsilon
$$
and
$$
\sup_{j\geq 1} |f(t-t_{j} )g(t-t_{j-1}) -  f(t-t_{j} )g(t-t_j)| \leq \epsilon .
$$
Hence Lemma \ref{mt5} shows that
\begin{eqnarray*}
&& \\ \lambda_t(f)\lambda_t(g)
&& \sim_{2\epsilon} \sum_{j=0}^{\infty} \Delta_j^4 f(t-t_j)g(t-t_j) +  \sum_{j=0}^{\infty} \Delta_j^2\Delta_{j+1}^2f(t-t_j)g(t-t_{j}) \\
&& \ \ \ \ \ \ \ \ \ \ \ \ \ \ \ \ \ \ \ \ \ \ \ \ \ +
\sum_{j=1}^{\infty} \Delta_{j}^2\Delta_{j-1}^2 f(t-t_{j} )g(t-t_{j}) \\
&&
= \lambda_t(fg)
\end{eqnarray*}
for all large enough $t$, i.e. (\ref{20}) holds. Note that $\lambda_t \circ \tau_s = \lambda_{t-s}$. Furthermore, it is obvious that $\lim_{t \to -\infty} \lambda_t(f) = 0$ for all $f \in S = C_0(\mathbb R)$, and that $\lim_{t \to \pm \infty} \lambda_t(g)b = g(\pm \infty)b$ for all $b \in B$ and all $g \in I$. In fact, $\lim_{t \to -\infty} \lambda_t(h) = h(-\infty)$, in norm. Furthermore, $\lambda_t(f) - f(-\infty)1 \in B$ for all $t,f$. To sum up:

\begin{lem}\label{a7} There is a completely positive and translation invariant asymptotic homomorphism $\lambda = \left(\lambda_t\right)_{t \in \mathbb R}: I \to M(B)$ such that
\begin{enumerate}
\item[1)] $\lambda_t(h) =\sum_{j=0}^{\infty} \Delta_j^2 h(t-t_j), \ t \in \mathbb R, h \in I$,
\item[2)] $\lambda_t(h) - h(-\infty)1 \in B , \ t \in \mathbb R, h \in I$,
\item[3)] $\lambda_t \circ \tau_s = \lambda_{t-s}, \ t,s \in \mathbb R$,
\item[4)] $\lim_{t \to -\infty} \lambda_t(h) = h(-\infty)1, \ h \in I$,
\item[5)] $\lim_{t \to \infty}\lambda_t(h)b = h(\infty)b, \ h \in I, b \in B$.
\end{enumerate}
\end{lem}

Since $A$ is separable there is a compact subset $X \subseteq A$ with dense span in $A$.

\begin{lem}\label{a5}

The sequence $h_0 \leq h_1 \leq h_2 \leq \cdots $ can be chosen so that
\begin{equation}\label{e6}
\sum_{n=0}^{\infty} \left\|\left[\Delta_n,s(a)\right]\right\| < \infty
\end{equation}
for all $a \in X$, and
\begin{equation}\label{e7}
s(a) - \sum_{j=0}^{\infty} \Delta_j s(a)\Delta_j \in B
\end{equation}
for all $a \in A$.
\end{lem}
\begin{proof} As is well-known there is a sequence $\delta_1 > \delta_2 > \delta_3 > \dots $in $]0,1[$ such that
$$
b \in B , \ 0 \leq b \leq 1, \  \|b s(a) - s(a)b\| < \delta_n \ \forall a \in X \  \Rightarrow \ \|\sqrt{b} s(a) - s(a)\sqrt{b}\| < 2^{-n} \ \forall a \in X .
$$
If we therefore choose $h_n$, as we can, such that $ \|u_n s(a) - s(a)u_n\| < {\delta_n}$ for all $n \in \mathbb N$ and all $a \in X$, we find that (\ref{e6}) holds. Note that
$$
s(a) -  \sum_{j=0}^{\infty} \Delta_j s(a)\Delta_j = \sum_{j=0}^{\infty} \Delta_j^2s(a) -  \Delta_j s(a)\Delta_j = \sum_{j=0}^{\infty} \Delta_j \left[\Delta_j, s(a)\right],
$$
for all $a \in A$, with strictly convergent sums. It follows from (\ref{e6}) that the last sum is convergent in norm when $a \in X$, so that $s(a) - \sum_{j=0}^{\infty} \Delta_j s(a)\Delta_j \in B$ for all $a \in X$. Since $s$ is linear modulo $B$ and since $X$ spans a dense subspace of $A$, we conclude that $s(a) - \sum_{j=0}^{\infty} \Delta_j s(a)\Delta_j \in B$ for all $a \in A$.
\end{proof}

Let $\kappa : \mathbb R \to ]0,1[$ be a continuous increasing function such that $\lim_{t \to - \infty} \kappa(t) = 0$ and $\lim_{t \to \infty} \kappa(t) = 1$. Set
$$
v_t = \sum_{j=0}^{\infty} \Delta_j^2 \kappa(t - t_j).
$$

\begin{lem}\label{a3} Assume that $h_0 \leq h_1 \leq h_2 \leq \dots$ is chosen so that (\ref{e6}) and (\ref{e7}) hold. Then

\begin{enumerate}
\item[a)] $\lim_{t \to -\infty}v_t = 0$,
\item[b)] $\lim_{t \to \infty} v_tb = b$ for all $b \in B$,
\item[c)] $\lim_{t \to \infty} v_ts(a) -s(a)v_t = 0$ for all $a\in A$,
\item[d)] $\lim_{t\to \infty}  g(v_t)s(a) - \varphi_t((g \circ \kappa) \otimes a) = 0$ for all $g \in C_0]0,1[, a \in A$.
\end{enumerate}
\end{lem}
\begin{proof} a) follows from Lemma \ref{mt5} because $\lim_{x \to -\infty} \kappa(x) = 0$. b) : Let $b \in B$ and $\epsilon > 0$ be given. There is an $N$ so large that $u_Nb \sim_{\epsilon} b$. Since
$$
\lim_{t \to \infty} v_tu_Nb = \lim_{t \to \infty} \sum_{j=0}^{N+1}\Delta_j^2\kappa(t-t_j)u_Nb =\sum_{j=0}^{N+1}\Delta_j^2u_Nb = u_Nb,
$$
we conclude that $v_tb \sim_{3 \epsilon} b$ for all $t$ large enough.

c): Since $s(a) - \sum_{j=0}^{\infty}\Delta_js(a)\Delta_j \in B$ for all $a \in A$, it suffices, thanks to b), to show that
\begin{equation}\label{e1}
\lim_{t \to \infty} v_t \left(\sum_{j=0}^{\infty}\Delta_js(a)\Delta_j\right) - \left(\sum_{j=0}^{\infty}\Delta_js(a)\Delta_j\right)v_t = 0.
\end{equation}
Let $\epsilon > 0$. Choose a $\delta > 0$ so small that $x,y \in \mathbb R, |x-y| < \delta \Rightarrow |\kappa(x) -\kappa(y)| < \epsilon$, and let $K \in \mathbb N$ be so large that $t_{n+2} - t_n < \delta$ when $n \geq K-1$. Set $\tilde{\kappa} = 1 - \kappa$ and note that $1 - v_t = \sum_{j=0}^{\infty} \Delta_j \tilde{\kappa}\left(t-t_j\right)\Delta_j$. Since $\lim_{x \to \infty} \tilde{\kappa}(x) = 0$ there is a $T \in \mathbb R$ so large that $\left\|1- v_t - \sum_{j=K}^{\infty} \Delta_j^2\tilde{\kappa}(t-t_j)\right\| < \epsilon$ when $t \geq T$. It follows that
\begin{equation*}
\begin{split}
&\left( 1 -v_t\right)\left(\sum_{j=0}^{\infty}\Delta_js(a)\Delta_j\right) \sim_{L \|a\| \epsilon} \sum_{j=K}^{\infty} \left( \Delta_j^2 \tilde{\kappa}(t - t_j) + \Delta_{j-1}^2\tilde{\kappa}(t-t_{j-1}) + \Delta_{j+1}^2\tilde{\kappa}(t - t_{j+1})\right) \Delta_js(a)\Delta_j \\
&\sim_{3 L\|a\| \epsilon} \sum_{j=K}^{\infty} \tilde{\kappa}(t-t_j)\Delta_js(a)\Delta_j
\end{split}
\end{equation*}
for all large enough $t$. Similarly, we find that
$$
\left(\sum_{j=0}^{\infty}\Delta_js(a)\Delta_j\right)\left(1 -v_t\right) \sim_{4L\|a\|\epsilon} \sum_{j=K}^{\infty} \tilde{\kappa}(t-t_j)\Delta_js(a)\Delta_j
$$
for all large enough $t$. Hence (\ref{e1}) holds.

d) : We may assume that $a \in X$. Since $s$ is homogeneous, we find that
$$
\varphi_t((g \circ \kappa)\otimes a) - \left(\sum_{j=0}^{\infty} \Delta_j^2g \circ \kappa(t-t_j)\right) s(a) = \sum_{j=0}^{\infty} g \circ \kappa(t - t_j)\Delta_j \left[s(a), \Delta_j\right] .
$$
Since $\sum_{j=0}^{\infty} \left\| \left[s(a), \Delta_j\right]\right\| < \infty$, it follows easily that
$$
\lim_{t \to \infty} \left[\varphi_t((g \circ \kappa)\otimes a) - \left(\sum_{j=0}^{\infty} \Delta_j^2g \circ \kappa(t-t_j)\right) s(a)\right] = 0.
$$
By Lemma \ref{a7}, $\lambda_t(h) = \sum_{j=0}^{\infty} \Delta_j^2h \circ \kappa(t-t_j)$ gives us an asymptotic homomorphism $\lambda : C_0(0,1] \to B$. It follows that
$$
\{h \in C_0(0,1]: \ \lim_{t \to \infty} h(v_t) - \lambda_t(h) = 0\}
$$
is a $C^*$-subalgebra of $C_0(0,1]$. Since it contains the identity function it must be all of $C_0(0,1]$, so we conclude that
$$
\lim_{t \to \infty} \left[g(v_t) s(a) - \left(\sum_{j=0}^{\infty} \Delta_j^2g \circ \kappa(t-t_j)\right) s(a)\right] = 0.
$$
\end{proof}

In conclusion: Let $CH(\psi) : SA \to B$ be the asymptotic homomorphism arising from the Connes-Higson construction introduced in \cite{CH}. Then $CH(\psi)$ is homotopic to a translation invariant asymptotic homomorphism.

\section{From translation invariant asymptotic homomorphisms to extensions.}\label{Sec I}

In this section we construct an inverse of $CH_{\tau}$.

Let $\beta_0 \in C_0(\mathbb R)$ be a continuous function with support in $[-1,1]$ such that
\begin{equation}\label{e12}
\beta_0(t) > 0, \ t \in ]-1,1[,
\end{equation}\begin{equation}
\tau_1(\beta_0)\tau_{-1}(\beta_0) = 0
\end{equation}
and
\begin{equation}\label{e11}
\tau_{-1}(\beta_0)\beta_0 + \beta_0^2 + \tau_{1}(\beta_0)\beta_0 = \beta_0.
\end{equation}
One choice of $\beta_0$ could be
$$
\beta_0(t) = \begin{cases} 0, & |t| \geq 1,\\ 1 + t, & t \in ]-1,0],\\ 1-t, & t \in [0,1[ . \end{cases}
$$
When $i \in \mathbb Z, \ i \neq 0$, set $\beta_i = \tau_{-i}\left(\beta_0\right)$. Then $\left\{\beta_i\right\}_{i \in \mathbb Z}$ is a partition of unity in $C_0(\mathbb R)$ such that $\beta_i$ is supported in $[-i-1,-i+1]$. Set $\alpha_i = \sqrt{\beta_i}$, and note that
\begin{equation}\label{17}
\alpha_i\alpha_j = 0, \ |i-j| \geq 2.
\end{equation}

Let $\varphi = \left(\varphi_t\right)_{t \in \mathbb R} :  SA \to B$ be an asymptotically translation invariant asymptotic homomorphism. For convenience we assume that $\varphi_t(0) = 0$ for all $t$. Then
\begin{equation}\label{ZZ2}
\varphi_i\left(\tau_{i}(\alpha_i \alpha_j)\otimes a\right) = 0 \ \text{when $|i-j| \geq 2$},
\end{equation}
for all $a \in A$. Let $\mathbb K$ be the $C^*$-algebra of compact operators on $l^2(\mathbb Z)$, and let $e_{ij}, i,j \in \mathbb Z$, be the standard matrix units in $\mathbb K$.

Define $\Phi_0: A \to M(B \otimes \mathbb K)$ by
$$
\Phi_0(a) = \sum_{i,j \in \mathbb Z} \varphi_i\left( \tau_{i}\left(\alpha_i \alpha_j\right)\otimes a\right) \otimes e_{ij}.
$$
The sum converges in the strict topology because of (\ref{ZZ2}) and $\sup_{i,j} \left\|\varphi_i\left( \tau_{i}\left(\alpha_i \alpha_j\right)\otimes a\right)\right\| < \infty$. Let $q_{B \otimes \mathbb K} : M(B \otimes \mathbb K) \to Q(B \otimes \mathbb K)$ be the quotient map, and set $\Phi = q_{B \otimes \mathbb K} \circ \Phi_0$.

\begin{lem}\label{l6} $\Phi$ is an extension of $A$ by $B \otimes \mathbb K$, i.e. $\Phi \in \Hom(A,Q(B \otimes \mathbb K))$.
\end{lem}
\begin{proof} Let $a,b \in A$. We check that $\Phi(a)\Phi(b) = \Phi(ab)$. Using (\ref{ZZ2}) we find that
\begin{equation}\label{e9}
\left\|\Phi(a)\Phi(b) - \Phi(ab)\right\|
\end{equation}
$$
\leq 3 \lim_{K \to \infty} \sup_{|i| \geq K,|k -i| \leq 3} \left\|\varphi_i\left(\tau_{i}(\alpha_i \alpha_k)\otimes ab\right) - \sum_{j \in \mathbb Z} \varphi_i\left( \tau_{i}(\alpha_i\alpha_j)\otimes a \right)  \varphi_j\left( \tau_{j}(\alpha_j\alpha_k)\otimes b \right) \right\| .
$$
Using (\ref{ZZ2}) again we find that
\begin{equation*}
\begin{split}
&\sum_{j \in \mathbb Z} \varphi_i\left(\tau_{i}(\alpha_i \alpha_j)\otimes a\right)\varphi_j\left(\tau_{j}(\alpha_j \alpha_k)\otimes b\right) \\
&=  \varphi_i\left(\tau_{i}(\alpha_i^2)\otimes a\right) \varphi_i\left(\tau_{i}(\alpha_i \alpha_k)\otimes b\right) + \varphi_i\left(\tau_{i}(\alpha_i\alpha_{i-1})\otimes a\right) \varphi_{i-1}\left(\tau_{i-1}(\alpha_{i-1} \alpha_k)\otimes b\right) \\
&+ \varphi_i\left(\tau_{i}(\alpha_i\alpha_{i+1})\otimes a\right) \varphi_{i+1}\left(\tau_{i+1}(\alpha_{i+1} \alpha_k)\otimes b\right)\\
& = \varphi_i\left(\alpha_0^2\otimes a\right) \varphi_i\left(\alpha_0 \alpha_{k-i}\otimes b\right) + \varphi_i\left(\alpha_0\alpha_{-1}\otimes a\right) \varphi_{i-1}\left(\alpha_{0} \alpha_{k-i+1}\otimes b\right) \\
&+ \varphi_i\left(\alpha_0\alpha_{1}\otimes a\right) \varphi_{i+1}\left(\alpha_{0} \alpha_{k-i-1}\otimes b\right).
\end{split}
\end{equation*}
Since $\varphi$ is translation invariant, we have that
$$
\lim_{i \to \infty} \sup_{|k-i| \leq 3} \left\|\varphi_i\left(\alpha_0\alpha_{-1}\otimes a\right) \varphi_{i-1}\left(\alpha_{0} \alpha_{k-i+1}\otimes b\right) - \varphi_i\left(\alpha_0\alpha_{-1}\otimes a\right) \varphi_{i}\left(\alpha_{-1} \alpha_{k-i}\otimes b\right)\right\|  = 0
$$
and
$$
\lim_{i \to \infty} \sup_{|k-i| \leq 3}\left\| \varphi_i\left(\alpha_0\alpha_{1}\otimes a\right) \varphi_{i+1}\left(\alpha_{0} \alpha_{k-i-1}\otimes b\right) -  \varphi_i\left(\alpha_0\alpha_{1}\otimes a\right) \varphi_{i}\left(\alpha_{1} \alpha_{k-i}\otimes b\right)\right\| = 0 .
$$

It follows that
\begin{eqnarray*}
&&\lim_{i \to \infty} \sup_{|k-i|\leq 3} \Bigl\|\sum_{j \in \mathbb Z} \varphi_i\left(\tau_{i}(\alpha_i \alpha_j)\otimes a\right)\varphi_j\left(\tau_{j}(\alpha_j \alpha_k)\otimes b\right)\Bigr. \\
&&\ \ \ \ \ \ \ \ \ \ \ \ \ \ \Bigl.- \varphi_i\left(\left(\alpha_0\alpha_0^2\alpha_{k-i} + \alpha_0 \alpha_{-1}^2\alpha_{k-i}+ \alpha_0 \alpha_{1}^2\alpha_{k-i}\right)\otimes ab \right)\Bigr\| = 0.
\end{eqnarray*}
Since $\alpha_0\alpha_0^2\alpha_{k-i} + \alpha_0 \alpha_{-1}^2\alpha_{k-i}+ \alpha_0 \alpha_{1}^2\alpha_{k-i} = \tau_{i}\left(\alpha_i\alpha_k\right)$ by (\ref{e11}), we see that
\begin{equation}\label{e8}
\lim_{i \to \infty} \sup_{|k-i|\leq 3} \left\|\sum_{j \in \mathbb Z} \varphi_i\left(\tau_{i}(\alpha_i \alpha_j)\otimes a\right)\varphi_j\left(\tau_{j}(\alpha_j \alpha_k)\otimes b\right) - \varphi_i\left(\tau_{i}(\alpha_i \alpha_k)\otimes ab\right)\right\| = 0.
\end{equation}
By using that $\lim_{t \to -\infty} \varphi_t(f) = 0$ for all $f \in SA$, we find easily that
\begin{equation}\label{e10}
\begin{split}
&
\lim_{i \to -\infty} \sup_{|k-i| \leq 3} \left\|\sum_{j \in \mathbb Z} \varphi_i\left(\tau_{i}(\alpha_i \alpha_j)\otimes a\right)\varphi_j\left(\tau_{j}(\alpha_j \alpha_k)\otimes b\right)\right\| \\
& = \lim_{i \to -\infty} \sup_{|k-i| \leq 3} \left\|\varphi_i\left(\tau_{i}(\alpha_i \alpha_k)\otimes ab\right)\right\| = 0.
\end{split}
\end{equation}
By combining (\ref{e9}), (\ref{e8}) and (\ref{e10}), we conclude that $\Phi(a)\Phi(b) = \Phi(ab)$.  It can be shown in a similar way that $\Phi$ is linear and self-adjoint. Thus $\Phi$ is an extension of $A$ by $B \otimes \mathbb K$.
\end{proof}

The extension $\Phi$ of Lemma \ref{l6} will be denoted by $I(\varphi)$. It is clear that the construction gives rise to a well-defined map
$$
I : [[SA,B]]_{a, \tau} \to \Ext_h(A,B \otimes \mathbb K)
$$
such that $I[\varphi] = [I(\varphi)]$.

\begin{rmk}\label{remork} \emph{At first sight it may seem odd that $I(\varphi)$ only depends of the 'discrete part' of the asymptotic homomorphism $\varphi$; i.e. of $\varphi_z, z \in \mathbb Z$. This is not so strange because a discrete translation invariant asymptotic homomorphism $\psi$, as defined in the obvious way, determines a unique element in $[[SA,B]]_{a,\tau}$ via the formula
$$
\psi_t = \psi_{[t]} \circ \tau_{[t] -t} .
$$}    
\end{rmk}

\begin{rmk}\label{remark} \emph{We want to point out that the $I$ construction, considered as a map $[[SA,B]]_{a, \tau} \to \Ext_h(A,B \otimes \mathbb K)$, can be wrapped in other, maybe more familiar guises: As pointed out in Remark \ref{remork} there is nothing lost in thinking of $\varphi$ as a translation invariant \emph{discrete} asymptotic homomorphism, and then the discrete version of the $M$-construction from \cite{H-LT} gives us a $\mathbb Z$-extension $E$ of $SA$ by $C_0(\mathbb Z,B)$. The crossed product $E \times \mathbb Z$ is then an extension of  $C(\mathbb T)\otimes \mathbb K  \otimes A$ by $B \otimes \mathbb K$. The pull back of this extension along the homomorphism $a \mapsto e \otimes a$, where $e \in C(\mathbb T)\otimes \mathbb K$ is a projection which generates $K_0( C(\mathbb T)\otimes \mathbb K) \simeq \mathbb Z$, gives us an  extension of $A$ by $B \otimes \mathbb K$ which is homotopic to $I(\varphi)$. Alternatively, the $I$ map can be realized via the genuine M-construction, using crossed products by $\mathbb R$.} 
\end{rmk}

\begin{lem}\label{l10} Let $s : B \to B \otimes \mathbb K$ be the stabilizing $*$-homomorphism, $s(b) = b \otimes e_{00}$, and let $\hat{s} : Q(B) \to Q(B \otimes \mathbb K)$ be the embedding induced by $s$. Define $s_* : \Ext_h(A,B) \to  \Ext_h(A, B \otimes \mathbb K)$ such that $s_*[\varphi] = [\hat{s} \circ \varphi]$. Then the diagram
\begin{equation*}
\begin{xymatrix}{
\Ext_h(A,B) \ar[r]^-{s_*} \ar[d]_-{CH_{\tau}} & \Ext_h(A, B \otimes \mathbb K) \\
[[SA,B]]_{a, \tau} \ar[ur]_-{I} & }
\end{xymatrix}
\end{equation*}
commutes.
\end{lem}
\begin{proof} We adopt the notation from above and from Section \ref{vch}. Let $\varphi : A \to Q(B)$ be an extension. Then $I \circ CH_{\tau}[\varphi]$ is represented by the map
$$
\Phi(a) = q_{B \otimes \mathbb K}\left( \sum_{i,j \in \mathbb Z} \sum_{k =0}^{\infty} \Delta_k^2 \alpha_i\alpha_j(-t_k) s(a) \otimes e_{ij} \right) .
$$
We may assume that there is a compact subset $X \subseteq A$ with dense span in $A$ such that
\begin{equation}\label{e15}
\sum_{k=0}^{\infty} \left\|\left[ \Delta_j, s(a)\right]\right\| < \infty
\end{equation}
for all $a \in X$, cf. Lemma \ref{a5}. Let $\lambda : I \to B$ be the completely positive asymptotic homomorphism $\lambda_t(h) = \sum_{j=0}^{\infty} \Delta_j^2 h(t-t_j)$, cf. Lemma \ref{a7}. Then
$$
 q_{B \otimes \mathbb K}\left( \sum_{i,j \in \mathbb Z} \sum_{k =0}^{\infty} \Delta_k^2 \alpha_i\alpha_j(-t_k) s(a) \otimes e_{ij} \right) =  q_{B \otimes \mathbb K} \left( \sum_{i,j \in \mathbb Z} \lambda_0(\alpha_i\alpha_j) s(a) \otimes e_{ij} \right) .
$$
Since $\lim_{k \to \infty} t_k - t_{k+1} = 0$, there is an $N \in \mathbb N$ so large that $\left|t_k - t_{k+1}\right| \leq N$ for all $k$. Since $\alpha_0$ is supported in $[-1,1]$, it follows from (\ref{e13}) that
$$
\lambda_i(\alpha_0)\lambda_j(\alpha_0) = \sum_{k,l =0}^{\infty} \Delta_k^2\Delta_l^2 \alpha_0(i-t_k)\alpha_0(j-t_l) = 0
$$
when $|i -j| \geq N+3$. Hence
$$
\sup_j \left\|\lambda_0(\alpha_i\alpha_j) - \lambda_0(\alpha_i)\lambda_0(\alpha_j)\right\| = \sup_j \left\|\lambda_i(\alpha_0\alpha_{j-i}) - \lambda_i(\alpha_0)\lambda_i(\alpha_{j-i})\right\|
$$
tends to $0$ as $i$ tends to $\pm \infty$. It follows that
$$
 q_{B \otimes \mathbb K}\left( \sum_{i,j \in \mathbb Z} \lambda_0(\alpha_i\alpha_j) s(a) \otimes e_{ij} \right)
=  q_{B \otimes \mathbb K}\left( \sum_{i,j \in \mathbb Z} \lambda_i(\alpha_0)\lambda_j(\alpha_0) s(a) \otimes e_{ij} \right).
$$
Set
$$
V = \sum_{i \in \mathbb Z} \lambda_0(\alpha_i) \otimes e_{0i} \in M(B \otimes \mathbb K),
$$
where the sum converges in the strict topology. By using that $\sum_{i \in \mathbb Z}\lambda_i(\alpha_0)^2 = 1$, with convergence in the strict topology, we see that $VV^* = 1 \otimes e_{00} \in M(B \otimes \mathbb K)$. Note that
$$
V^*\left(s(a) \otimes e_{00} \right) V =  \sum_{i,j \in \mathbb Z} \lambda_i(\alpha_0)s(a)\lambda_j(\alpha_0)  \otimes e_{ij} .
$$
It follows from (\ref{e15}) that $\lim_{i \to \infty} \lambda_i(\alpha_0)s(a) - s(a)\lambda_i(\alpha_0) = 0$ for all $a \in X$. Since $\lim_{i \to -\infty} \lambda_i(\alpha_0) = 0$, we conclude that
$$
V^* \left(s(a) \otimes e_{00} \right)V  =   \sum_{i,j \in \mathbb Z} \lambda_i(\alpha_0)\lambda_j(\alpha_0) s(a) \otimes e_{ij} 
$$
modulo $B \otimes \mathbb K$ for all $a \in X$. Then $U = \left( \begin{smallmatrix} V & 1-VV^*\\ 1-V^*V & V^* \end{smallmatrix} \right) \in M_2(M(B \otimes \mathbb K))$ is a unitary such that
\begin{equation}\label{e16}
q_{M_2(B \otimes \mathbb K)} \left( \Ad U \circ  \left( \begin{smallmatrix} \Phi(a) & \\  & 0 \end{smallmatrix} \right)\right)  =q_{M_2(B \otimes \mathbb K)}  \left( \begin{smallmatrix} s(a) \otimes e_{00} &  \\  & 0 \end{smallmatrix} \right)
\end{equation}
for all $a \in X$. Since $X$ spans a dense set in $A$, (\ref{e16}) holds for all $a \in A$. Now the lemma follows because $U$ can be connected to $1$ in the strict topology, by a path of unitaries in $M(M_2(B \otimes \mathbb K))$.
\end{proof}

\begin{lem}\label{Lem2} Let $\varphi = \left( \varphi_t\right)_{t \in \mathbb R} : SA \to B$ be an asymptotically translation invariant asymptotic homomorphism. It follows that there is a homogeneous equi-continuous and uniformly continuous asymptotic homomorphism $\psi =  (\psi_t)_{t \in \mathbb R} : SA \to B$ such that
\begin{enumerate}
\item[1)] $\lim_{t \to \infty} \psi_t(f) - \varphi_t(f) = 0$ for all $f \in SA$,
\item[2)] $\lim_{t \to \infty} \sup_{s \in K} \left\|\psi_{t-s}(f) - \psi_t\left(\tau_s(f)\right)\right\| = 0$ for all $f \in SA$ and all compact subsets $K \subseteq \mathbb R$,
\item[3)] $\lim_{t \to -\infty} \psi_t(f) = 0$ for all $f \in SA$.
\end{enumerate}

\end{lem}
\begin{proof} Let $X = \left\{ f \in C_b(\mathbb R,B) : \ \lim_{t \to - \infty} f(t) = 0\right\}$ and $X_0 = \left\{f \in X: \ \lim_{t \to \infty} f(t) = 0\right\}$. Define an automorphism $\theta_x$ of $X$ such that $\theta_x(f)(t) = f(t-x)$, and note that $\theta_x$ leaves $X_0$ globally invariant, so that $\theta_x$ induces an automorphism of $X/X_0$ which we denote by $\theta_x$ again. Let $\Phi_0 : SA \to X$ be the map $\Phi_0(f)(t) = \varphi_t(f)$. Let $q : X \to X/X_0$ 
be the quotient map, and note that $\Phi = q \circ \Phi_0$ is then an equivariant $*$-homomorphism. Let $X_{\mathbb R}$ denote the $C^*$-subalgebra of $X$ consisting of the elements $f \in X$ for which $x \mapsto \theta_x(f)$ is continuous. It follows from Theorem 2.1 of \cite{Th?} that $ \Phi(SA) \subseteq q\left(X_{\mathbb R}\right)$. Let $L : q\left( X_{\mathbb R}\right)  \to X_{\mathbb R}$ be a continuous and 
homogeneous right-inverse for $q : X_{\mathbb R} \to q\left( X_{\mathbb R}\right)$, and set $\psi_t(f) = L \circ \Phi(f)(t)$. Then $\psi = \left( \psi_t\right)_{t \in \mathbb R}$ is an equicontinuous and homogeneous asymptotic homomorphism such that $\lim_{t \to \infty} \psi_t(f) - \varphi_t(f) = 0$, $\lim_{t \to - \infty} \psi_t(f) = 0$ and
\begin{equation}\label{16}
\lim_{t \to \infty} \psi_t\left(\tau_s(f)\right) - \psi_{t-s}(f) = 0
\end{equation}
for all $f \in SA$ and all $s \in \mathbb R$. The uniform continuity of $\psi$ follows because $L \circ \Phi(f) \in X_{\mathbb R}$. And then 2) follows from (\ref{16}), the uniform continuity and the equicontinuity of $\psi$, as follows: By uniform continuity there is $\delta > 0$ so small that $\sup_{t \in \mathbb R}\|\psi_{t-x}(f) - \psi_{t-y}(f)\| \leq \epsilon$ when $|x-y| < \delta$. By equicontinuity of $\psi$ there is for each $s \in K$, a $\delta_s \in ] 0, \delta[$ such 
that $\sup_{t \in \mathbb R} \left\|\psi_{t}(\tau_s(f)) - \psi_t(\tau_{s'}(f))\right\| < \epsilon$ when $|s'-s| < \delta_s$. Let $\left]s_i - \delta_{s_i}, s_i + \delta_{s_i}\right[, i = 1,2, \dots, L$, be a finite subcover of the cover $]s - \delta_s,s + \delta_s[, s \in K$, of $K$. Choose $T > 0$ so large that $\left\|\psi_{t-s_i}(f) - \psi_t \circ \tau_{s_i}(f)\right\| < \epsilon$ for all $i = 1,2, \dots, L$, when $t > T$. It follows then that
$$
 \sup_{s \in K} \left\|\psi_{t-s}(f) - \psi_t\left(\tau_s(f)\right)\right\| < 3 \epsilon
$$
when $t > T$.
\end{proof}

It follows from Lemma \ref{Lem2} that we can always replace an asymptotically translation invariant asymptotic homomorphism with one which has the properties spelled out in Lemma \ref{Lem2}; namely, uniform continuity and equicontinuity. Property 2) of Lemma \ref{Lem2} is then automatic.

\begin{lem}\label{pre} Let $\varphi = \left( \varphi_t\right)_{t \in \mathbb R} : SA \to B$ be an equi-continuous asymptotically translation invariant asymptotic homomorphism, and let $M \subseteq SA$ be a pre-compact subset. It follows that
\begin{enumerate}
\item[i)] $\lim_{t \to - \infty} \sup_{x \in M} \left\|\varphi_t(x)\right\| = 0$.
\item[ii)] For all $\epsilon > 0$, there is a $T > 0$ such that
$$
\left\|\varphi_t(x) -\varphi_t(y)\right\| \leq \| x-y\| + \epsilon
$$
for all $x,y \in M$ and all $t \geq T$.
\item[iii)]  For all $\epsilon > 0$, there is a $T > 0$ such that
$$
\left\|\varphi_t(x)\varphi_t(y) - \varphi_t(xy)\right\| \leq \epsilon
$$
for all $x,y \in M$ and all $t \geq T$.
\item[iv)]  For all $\epsilon > 0$, there is a $T > 0$ such that
$$
\left\|\varphi_t(x) + \varphi_t(y) - \varphi_t(x+y)\right\| \leq \epsilon
$$
for all $x,y \in M$ and all $t \geq T$.
\item[v)]  For any $\epsilon > 0$, and any compact subset $K \subseteq \mathbb R$, there is a $T > 0$ such that
$$
\sup_{s \in K} \left\|\varphi_{t-s}(x) - \varphi_t(\tau_s(x)))\right\| \leq \epsilon
$$
for all $x\in M$ and all $t \geq T$.
\end{enumerate}
\end{lem}
\begin{proof} 
Left to the reader.
\end{proof}

\begin{lem}\label{l12}  The diagrams
\begin{equation*}
\begin{xymatrix}{
& \Ext_h(A, B \otimes \mathbb K) \ar[d]^-{CH_{\tau}} \\
[[SA,B]]_{a, \tau} \ar[ur]^-I  \ar[r]_-{s_*} & [[SA,B \otimes \mathbb K]]_{a, \tau} }
\end{xymatrix}
\end{equation*}
and
\begin{equation*}
\begin{xymatrix}{
& \Ext_h(A, B \otimes \mathbb K) \ar[d]^-{CH_{\tau}} \\
[[SA,B]]_{ \tau} \ar[ur]^-I  \ar[r]_-{s_*} & [[SA,B \otimes \mathbb K]]_{ \tau} }
\end{xymatrix}
\end{equation*}
both commute.
\end{lem}
\begin{proof}  As in the previous proof we will adopt the notation from above and from Section \ref{vch}. Let $\varphi = \left(\varphi_t\right)_{t \in \mathbb R} : SA \to B$ be an asymptotically translation invariant asymptotic homomorphism. By Lemma \ref{Lem2} we may assume that $\varphi$ has the properties spelled out there (for $\psi$). Set $\Lambda(a) = \sum_{i,j \in \mathbb Z} \varphi_i\left(\tau_i\left(\alpha_i\alpha_j\right)\otimes a\right) \otimes e_{ij}$. Let $c$ be a strictly positive element in $B$. Then
$$
b = \sum_{i \in \mathbb Z} 2^{-|i|}c \otimes e_{ii}
$$
is a strictly positive element of $B \otimes \mathbb K$ which we may use to construct $CH_{\tau}\left(I(\varphi)\right)$. However, we consider instead
$$
b_0 = \sum_{i \in \mathbb Z} 2^{-|i|} \otimes e_{ii},
$$
which in general is not an element of $B \otimes \mathbb K$, only of $M(B \otimes \mathbb K)$. But $b_0$ has the following two properties: $b_0\Phi_0(A) \subseteq B \otimes \mathbb K$ and $\lim_{k \to \infty} b_0^{\frac{1}{k}}x = x$ for all $x \in B \otimes \mathbb K$. We can therefore find a sequence $h_0 \leq h_1 \leq h_2 \leq \dots $ of functions in $C_0]0,1]$ such that 
\begin{equation}\label{u1}
h_{n+1}h_n = h_n
\end{equation} 
for all $n$, and the sequence $u_n = h_n(b_0), n = 0,1,2, \dots$ has the properties 
\begin{enumerate}
 \item[a)] $\lim_{n \to \infty} u_n\Phi_0(a) - \Phi_0(a)u_n =0$ for all $a \in A$, and \item[b)] $\lim_{n \to \infty} u_nx = x$ for all $x \in B \otimes \mathbb K$. 
\end{enumerate}
Set $\Delta_0 = \sqrt{u_0}, \Delta_n = \sqrt{u_n - u_{n-1}}, n \geq 1$. The arguments from Section \ref{subsec} show that $CH_{\tau} \circ I[\varphi]$ is represented by a translation invariant asymptotic homomorphism $\psi : SA \to B \otimes \mathbb K$ such that
\begin{equation}\label{u3}
\lim_{t \to \infty} \left[\psi_t(f \otimes a) - \Phi_0(a)\sum_{j=0}^{\infty} \Delta_j^2 f(t-t_j)\right] = 0 .
\end{equation}
Note that we can choose the sequence $\{t_j\}_{j=0}^{\infty}$ in $[0,\infty)$ at will, as long as $\lim_{j \to \infty} t_j = \infty$ and $\lim_{j \to \infty} t_{j+1} - t_j = 0$. This freedom is used as follows: For each $i \in 0,1,2,3, \dots $, set
$$
j_i = \min \left\{ k : h_k\left(2^{-i}\right) = 1 \right\} .
$$
It follows from (\ref{u1}) that $j_i \leq j_{i+1}$ and from b) that $\lim_{i \to \infty} j_i = \infty$. We can therefore choose $\{t_j\}_{j=0}^{\infty}$ such that $\lim_{i \to \infty} t_{j_{i}} - t_{j_{i-1}} = 0$. With this choice, it follows from (\ref{u1}) that
\begin{equation}\label{u2}
\begin{split}
&\sum_{j=0}^{\infty} \Delta_j^2 f(t-t_j) = \\
& \sum_{i \in \mathbb Z} \left(\left(1-h_{j_{|i|} -1}\left(2^{-|i|}\right)\right) f\left(t-t_{j_{|i|}}\right) +  h_{j_{|i|} -1}\left(2^{-|i|}\right)f\left(t - t_{j_{|i|} -1}    \right) \right) e_{ii} ,
\end{split}
\end{equation}
if we set $h_{-1} = 0$. Set $s_i = t_{j_{|i|}}, i \in \mathbb Z$, and note that $\lim_{i \to \pm \infty} \left|s_i - s_{i-1}\right| = 0$ and $\lim_{i \to \pm \infty} s_i = \infty$. It follows from (\ref{u2}) that
$$
\lim_{t \to  \infty} \left[\sum_{j=0}^{\infty} \Delta_j^2 f(t-t_j) - \sum_{i \in \mathbb Z} f(t-s_i)e_{ii}\right] = 0
$$
for all $f \in C_0(\mathbb R)$, and then from (\ref{u3}) that
$$
\lim_{t \to \infty} \left[\psi_t(f \otimes a) -  \sum_{i,j \in \mathbb Z} \varphi_i\left( \tau_{i}\left(\alpha_i \alpha_jf(t-s_i )\right)\otimes a\right) \otimes e_{ij}\right] = 0,
$$
for all $f \in C_0(\mathbb R), a \in A$. Let $\beta : \mathbb R \to \mathbb R$ be a continuously differentiable decreasing function such that $\beta(i) = s_i =t_{j_{|i|}}$ for all $i\in \mathbb Z$, $i < 0$, while $\beta(x) = -x, x \geq 0$, and such that $\beta^{\prime}$ is bounded. The last condition can be met because  $\lim_{i \to \pm \infty} \left|s_i - s_{i-1}\right| = 0$. Then
\begin{equation}\label{u5}
\lim_{t \to \infty} \left[\psi_t(f \otimes a) -  \sum_{i,j \in \mathbb Z} \varphi_i\left( \tau_{i}\left(\alpha_i \alpha_jf( t -\beta(\cdot))\right)\otimes a\right) \otimes e_{ij}\right] = 0.
\end{equation}
Indeed, we claim that
\begin{equation}\label{u4}
\lim_{t \to \infty} \sup_{i,j \in \mathbb Z} \left\| \varphi_i\left( \tau_{i}\left(\alpha_i \alpha_jf(t-s_i )\right)\otimes a\right) - \varphi_i\left( \tau_{i}\left(\alpha_i \alpha_jf( t -\beta(\cdot))\right)\otimes a\right)\right\| = 0.
\end{equation}
To see this, let $\epsilon > 0$ and observe that $\tau_{i}\left(\alpha_i \alpha_jf(t-s_i )\right) = \alpha_0\alpha_{j-i} f(t-s_i)$ is a pre-compact set in $C_0(\mathbb R)$. By i) of Lemma \ref{pre} there is a $K < 0$ such that 
\begin{equation}\label{u8}
 \left\| \varphi_i\left( \tau_{i}\left(\alpha_i \alpha_jf(t-s_i )\right)\otimes a\right)\right\| \leq \epsilon
\end{equation}
for all $i \leq K$. Similarly, it follows from Ascoli's theorem, cf. e.g. Theorem A 5 on page 369 of \cite{R}, that  
$$
\left\{ \tau_{i}\left(\alpha_i \alpha_jf( t -\beta(\cdot))\right) : i \in \mathbb Z, t \in \mathbb R \right\} = \left\{ \alpha_0 \alpha_{j-i}f( t -\beta(\cdot - i)) : i \in \mathbb Z, t \in \mathbb R \right\}
$$
is pre-compact in $C_0(\mathbb R)$. By using i) of Lemma \ref{pre} again, we may assume, by decreasing $K$ if necessary, that also  
\begin{equation}\label{u9}
\left\| \varphi_i\left( \tau_{i}\left(\alpha_i \alpha_jf(t- \beta(\cdot ) )\right)\otimes a\right)\right\| \leq \epsilon 
\end{equation}
 for all $i \leq K$. It follows from the uniform continuity of $f$, and the fact that $\lim_{i \to \infty} \sup_{x \in [-1,1]} |s_i -\beta(x - i)| = 0$, that
$$
\lim_{i \to \infty} \sup_{t,j} \left\| \alpha_0\alpha_{j-i}f(t - s_i) -  \alpha_0\alpha_{j-i}f(t - \beta(\cdot - i))\right\| = 0.
$$
It follows therefore from ii) of Lemma \ref{pre} that there is an $L > 0$ so large that
\begin{equation}\label{u6}
\left\| \varphi_i\left( \tau_{i}\left(\alpha_i \alpha_jf(t-s_i )\right)\otimes a\right) - \varphi_i\left( \tau_{i}\left(\alpha_i \alpha_jf( t -\beta(\cdot))\right)\otimes a\right)\right\| \leq \epsilon
\end{equation}
when $i \geq L$. Since 
$$
\lim_{t \to \infty} \sup_{K \leq i\leq L} \left\| \tau_{i}\left(\alpha_i \alpha_jf(t-s_i )\right)\right\| = \lim_{t \to \infty}  \sup_{K \leq i\leq L} \left\|\tau_{i}\left(\alpha_i \alpha_jf( t -\beta(\cdot))\right)\right\| = 0,
$$ 
we conclude from (\ref{u8}),(\ref{u9}) and (\ref{u6}) that (\ref{u4}) holds. Hence (\ref{u5}) holds.

To proceed it is useful to apply the following statement, which is an abstract version of Ascoli's theorem. In fact, the proof is identical to the standard proof of Ascoli's theorem, as presented for example in \cite{R}. Alternatively, it can be deduced from the version of Ascoli's theorem on page 81 of \cite{Ke}.

\begin{thm}\label{aA} Let $K$ be a locally compact Hausdorff space, $X$ a Banach space, and $C_0(K,X)$ the Banach space of continuous $X$-valued functions on $K$ that vanish at infinity, with the norm $\|f\| = \sup \left\{\|f(k)\|: \ k \in K\right\}$. A subset $M \subseteq C_0(K,X)$ is pre-compact (i.e. $\overline{M}$ is compact) in $C_0(K,X)$ if and only if
\begin{enumerate}
\item[a)] $\overline{ \{f(k): f \in M\}}$ is compact in $X$ for all $k \in X$,
\item[b)] $M$ is equi-continuous, i.e. for all $\epsilon > 0$, every $k \in K$ has a neighbourhood $V$ such that $\|f(k) - f(k')\| < \epsilon$ for all $k' \in V$ and all $f \in M$,
\item[c)] the elements of $M$ vanish equicontinuously at infinity, i.e. for all $\epsilon > 0$ there is compact subset $C \subseteq K$ such that $\|f(k)\| \leq \epsilon$ for all $f \in M$ and all $k \in K \backslash C$.
\end{enumerate}
\end{thm}

Set now $\beta_{\lambda}(x) = \lambda \beta(x) - (1-\lambda)x$, and note that $\sup_{\lambda \in [0,1]} \left\|\beta_{\lambda}^{\prime}\right\| < \infty$. It follows then from Theorem \ref{aA} that for any $f \in SA = C_0(\mathbb R, A)$, the set 
$$
\left\{ \alpha_0\alpha_{j-i}f\left(t - \beta_{\lambda}( \cdot - i)\right): \ j,i \in \mathbb Z, t \in \mathbb R, \lambda \in [0,1] \right\}
$$
is a pre-compact set in $SA$. We claim that as a consequence of this we have that
\begin{align}\label{uu1}
\lim_{t \to \infty} \sup_{i,k \in \mathbb Z} \sup_{\lambda \in [0,1]} \bigg\|\sum_{j \in
    \mathbb Z}  \varphi_i\left( \tau_{i}\left(\alpha_i
      \alpha_jf(t - \beta_{\lambda}(\cdot))\right)\right)  \varphi_j\left(
    \tau_{j}\left(\alpha_j \alpha_kg(t -\beta_{\lambda}(\cdot))\right)\right)
  \\-  \varphi_i\left( \tau_{i}\left(\alpha_i \alpha_kfg(t -\beta_{\lambda}(\cdot)
   )\right)\right)\bigg\| = 0
\end{align}
for all $f,g \in SA$.

To see this, note that by iii) and iv) of Lemma \ref{pre}, there is a $K > 0$ such that
\begin{align}\label{uu2}
  \bigg\|\sum_{j \in
    \mathbb Z}  \varphi_i\left( \tau_{i}\left(\alpha_i
      \alpha_jf(t - \beta_{\lambda}(\cdot))\right)\right)  \varphi_j\left(
    \tau_{j}\left(\alpha_j \alpha_kg(t -\beta_{\lambda}(\cdot))\right)\right)
  \\-  \varphi_i\left( \tau_{i}\left(\alpha_i \alpha_kfg(t -\beta_{\lambda}(\cdot)
   )\right)\right)\bigg\| \leq \epsilon
\end{align}
for all $\lambda \in [0,1],k \in \mathbb Z$ when $i \geq K$. Since $\epsilon > 0$ is arbitrary, (\ref{uu1}) will follow if we can show that there is a $T > 0$ so large that 
\begin{equation}\label{uu3}
 \left\|\varphi_i\left( \tau_{i}\left(\alpha_i
      \alpha_jf(t - \beta_{\lambda}(\cdot))\right)\right)\right\| \leq \epsilon
\end{equation}    
for all $i \leq K+1$ and all $\lambda \in [0,1]$ when $t \geq T$. Since 
\begin{equation}\label{uu14}
\left\{ \alpha_0 \alpha_{j-i} f(t - \beta_{\lambda} ( \cdot - i)) : \ i,j \in \mathbb Z, t \in \mathbb R, \lambda \in [0,1]\right\}
\end{equation}
is a pre-compact set in $SA$ by Theorem \ref{aA}, it follows from i) of Lemma \ref{pre} that there is a $L < 0$ such that (\ref{uu3}) holds for all $t \in \mathbb R$ and all $\lambda \in [0,1]$ when $i \leq L$. Since
$$
\lim_{t \to \infty} \sup \left\{ \left\| \alpha_0(s) \alpha_{j-i}(s) f(t - \beta_{\lambda} ( s  - i)) \right\|  : \ s \in [-1,1], L \leq i \leq K \right\} = 0,
$$
the equi-continuity of $\varphi$ at $0$ ensures that 
$$
\lim_{ t \to \infty} \sup \left\{\left\|\varphi_i\left( \tau_{i}\left(\alpha_i
      \alpha_jf(t - \beta_{\lambda}(\cdot))\right)\right)\right\|:  \ \lambda \in [0,1], \ L  \leq i \leq K \right\} = 0.
$$
We can therefore find $T > 0$ such that (\ref{uu3}) holds, proving (\ref{uu1}). Similar arguments show that 
\begin{align}\label{uu11}
\lim_{t \to \infty} \sup_{i,j \in \mathbb Z} \sup_{\lambda \in [0,1]} \bigg\|  \varphi_i\left( \tau_{i}\left(\alpha_i
      \alpha_jf(t - \beta_{\lambda}(\cdot))\right)\right)  + \mu \varphi_i\left(
    \tau_{i}\left(\alpha_i \alpha_jg(t -\beta_{\lambda}(\cdot))\right)\right)
  \\-  \varphi_i\left( \tau_{i}\left(\alpha_i \alpha_j\left(f + \mu g\right)(t -\beta_{\lambda}(\cdot)
   )\right)\right)\bigg\| = 0,
\end{align}
and that
\begin{equation}\label{uu12}
\lim_{t \to \infty} \sup_{i,j \in \mathbb Z} \sup_{\lambda \in [0,1]} \bigg\|  \varphi_i\left( \tau_{i}\left(\alpha_i
      \alpha_jf(t - \beta_{\lambda}(\cdot))^*\right)\right) -  \varphi_j\left( \tau_{j}\left(\alpha_i \alpha_jf(t -\beta_{\lambda}(\cdot)
   )\right)\right)^*\bigg\| = 0,
\end{equation}
for all $f,g \in SA, \mu \in \mathbb C$. For $t \in \mathbb R$, define $\Phi_t : SA \to IB$ by
$$
\Phi_t(f)(\lambda) =  \sum_{i,j \in \mathbb Z} \varphi_i\left( \tau_{i}\left(\alpha_i \alpha_jf( t -\beta_{\lambda}(\cdot))\right)\right) \otimes e_{ij} .
$$
The continuity in $\lambda$ follows from another application of Theorem \ref{aA} and ii) of Lemma \ref{pre}. To see that that $\lim_{t \to - \infty} \Phi_t(f) = 0$ use first the pre-compactness of (\ref{uu14}) to find a $K < 0$ such that 
\begin{equation}\label{ZZ3}
\sup_{j,t,\lambda} \left\|\varphi_i\left( \tau_i (\alpha_i \alpha_j f(t - \beta_{\lambda}( \cdot)) )\right)\right\| \leq \epsilon
\end{equation} 
when $i \leq K$. Then observe that there is a $c \in \mathbb R$ such that $\beta_{\lambda}(x - i) \geq c$ for all $\lambda \in [0,1]$, all $x \in [-1,1]$ and all $i \geq K$. It follows that 
$$
\lim_{t \to -\infty} \left[\sup_{\lambda, i \geq K} \left\| \alpha_0\alpha_{j-i}f\left(t - \beta_{\lambda}( \cdot - i)\right)\right\| \right] = 0.
$$ 
The equi-continuity of $\varphi$ at $0$ implies then that
$$
\lim_{t \to -\infty} \left[\sup_{\lambda, i \geq K} \left\| \varphi_i \left( \tau_i \left(\alpha_i\alpha_{j}f\left(t - \beta_{\lambda}( \cdot)\right)\right)\right)\right\| \right] = 0.
$$     
Combined with (\ref{ZZ3}) this shows that $\lim_{t \to - \infty} \Phi_t(f) = 0$, as desired. Together with (\ref{uu1}), (\ref{uu11}) and (\ref{uu12}) we see that $\left(\Phi_t\right)_{t \in \mathbb R} : SA \to IB$ is an asymptotic homomorphism. We claim that $\Phi$ is translation invariant. The algebraic condition, that $\Phi_{t-s} = \Phi_t \circ \tau_s$, is trivially satisfied, so it remains to show that $\Phi$ is equi-continuous. To prove this, let $C \subseteq SA$ be a compact subset containing the set (\ref{uu14}). Let $\epsilon > 0$. By equi-continuity of $\varphi$ and compactness of $C$ there is a finite set $g_1,g_2, \dots, g_N \in C$ and a $\delta > 0$ such that $\bigcup_i \{h \in SA: \ \|h - g_i\| < \frac{\delta}{2}\} \supseteq C$, and $\sup_t \|\varphi_t(h) - \varphi_t(g_i)\| \leq \epsilon$ when $\|h - g_i\| < \delta$. It follows $\sup_t \|\Phi_t(h) - \Phi_t(f)\| \leq  3\epsilon$ when $\|h - f\| \leq \frac{\delta}{2}$, proving that $\Phi$ is indeed equi-continuous. Thus $\Phi$ gives us a homotopy of translation invariant asymptotic homomorphism which thanks to (\ref{u5}) connects $\psi$ with the translation invariant asymptotic homomorphism $\psi' : SA \to B$ such that
$$
\psi'_t(f) = \sum_{i,j \in \mathbb Z} \varphi_i \left( \tau_i \left( \alpha_i \alpha_j f(t + \cdot)\right)\right) \otimes e_{ij},
$$
$f \in SA$. Note that 
$$
\sum_{i,j \in \mathbb Z} \varphi_i \left( \tau_i \left( \alpha_i \alpha_j f(t + \cdot)\right)\right) \otimes e_{ij} = \sum_{i,j \in \mathbb Z} \varphi_{t + (i-t)} \left(  \alpha_0 \alpha_{j-i} f(t + \cdot - i)\right) \otimes e_{ij}.
$$
Let $\epsilon > 0$. Since $\varphi$ is equi-continuous at $0$ there is $\delta > 0$ such that $\|x \| \leq \delta \Rightarrow \sup_t \left\|\varphi_t(x)\right\| \leq \epsilon$. Since $f$ vanishes at infinity there is a $K > 0$ such that $\|f(s)\| \leq \delta$ when $|s| \geq K-1$. It follows that $\left\| \alpha_0 \alpha_{j-i} f(t + \cdot - i)\right\| \leq \delta $ when $|t-i| \geq K$, and hence that 
$$
\left\|\sum_{i,j \in \mathbb Z} \varphi_i \left( \tau_i \left( \alpha_i \alpha_j f(t + \cdot)\right)\right) \otimes e_{ij} - \sum_{i,j \in \mathbb Z, |t-i| \leq K }  \varphi_{t + (i-t)} \left(  \alpha_0 \alpha_{j-i} f(t + \cdot - i)\right) \otimes e_{ij}\right\| \leq 3 \epsilon .
$$
It follows from Theorem \ref{aA} that
$$
\left\{ \alpha_0 \alpha_{j-i} f(t + \cdot - i) : \ t \in \mathbb R, i,j \in \mathbb Z\right\}
$$
is a pre-compact subset of $SA$, so we can apply v) Lemma \ref{pre} to conclude that
$$
\left\| \sum_{i,j \in \mathbb Z, |t-i| \leq K }  \varphi_{t + (i-t)} \left(  \alpha_0 \alpha_{j-i} f(t + \cdot - i)\right) \otimes e_{ij} -  \sum_{i,j \in \mathbb Z, |t-i| \leq K }  \varphi_{t} \left( \tau_t\left(  \alpha_i \alpha_{j} f(t + \cdot)\right)\right) \otimes e_{ij}\right\| \leq \epsilon
$$
for all $t$ large enough. Since $\left\| \tau_t\left(  \alpha_i \alpha_{j} f(t + \cdot)\right)\right\| = \left\|\alpha_0\alpha_{j-i} f(t + \cdot -i)\right\| \leq \delta$ when $|t-i|\geq K$, we find all together that
$$
\left\|\psi'_t(f) - \sum_{i,j \in \mathbb Z }  \varphi_{t} \left( \tau_t\left(  \alpha_i \alpha_{j} f(t + \cdot)\right)\right) \otimes e_{ij}\right\| \leq 7\epsilon
$$
for all $t $ large enough. It follows from Theorem \ref{aA} that
$$
\left\{ \tau_t(\alpha_i\alpha_j)f: \ t \in \mathbb R, i,j \in \mathbb Z\right\}
$$
is a pre-compact subset of $SA$, so we conclude from i) of Lemma \ref{pre} that 
$$
\lim_{t \to - \infty}  \sum_{i,j \in \mathbb Z }  \varphi_{t} \left( \tau_t\left(  \alpha_i \alpha_{j} f(t + \cdot)\right)\right) \otimes e_{ij} = 0 .
$$
It follows that if we set
$$
\psi''_t(f)  =  \sum_{i,j \in \mathbb Z }  \varphi_{t} \left( \tau_t\left(  \alpha_i \alpha_{j}\right) f\right) \otimes e_{ij} ,
$$
then $\psi''=\left(\psi''_t\right)_{t \in \mathbb R} : SA \to B$ is a translation invariant asymptotic homomorphism which asymptotically agrees with $\psi'$ and hence in particular is homotopic to $\psi'$.

Now define a continuous path $\{\beta^{\lambda}_i\}_{i \in \mathbb Z}, \lambda \in [0,\infty[$, of partitions of unity in $C_0(\mathbb R)$ such that
$$
\beta^{\lambda}_i(\cdot) = \alpha^2_i(\cdot - \lambda), \quad i \geq 1, \quad
\beta^{\lambda}_i(\cdot) = \alpha^2_i(\cdot + \lambda), \quad i \leq -1,
$$
and
$$
\beta^{\lambda}_0(\cdot) = 1 - \sum_{i=1}^{\infty} \left( \alpha_i^2(\cdot - \lambda) + \alpha_{-i}^2(\cdot + \lambda) \right) .
$$
Set
$$
\alpha_i^{\lambda} = \sqrt{\beta_i^{\frac{\lambda}{1-\lambda}}}, \lambda \in [0,1[,
$$
and 
$$
\alpha_i^1 = 0, i \neq 0, \ \  \alpha_0^1 = 1.
$$
Note that
$$
\alpha_i^{\lambda}\alpha_j^{\lambda} = 0, \ |i -j| \geq 2,
$$
and 
\begin{equation}\label{uu13}
\sum_{ j \in \mathbb Z} \left(\alpha_j^{\lambda}\right)^2 \alpha_k^{\lambda} = \alpha_k^{\lambda}
\end{equation}
for all $\lambda \in [0,1]$ and $k \in \mathbb Z$. For $t \in \mathbb R, \lambda \in [0,1]$, set
$$
\Psi_t(f)(\lambda) =  \sum_{i,j \in \mathbb Z} \varphi_t\left( \tau_{t} \left(\alpha_i^{\lambda} \alpha_j^{\lambda}\right)f\right) \otimes e_{ij} .
$$
Since $f$ vanishes at infinity we have that
$$
\lim_{i \to \pm \infty}  \sup \left\{ \left\|\tau_{t} \left(\alpha_i^{\lambda} \alpha_j^{\lambda}\right)f\right\|: \  j \in \mathbb Z, \lambda \in [0,1] \right\} = 0,
$$
so the continuity of $\varphi_t$ at $0$ ensures that for any $\epsilon > 0$, there is a $K > 0$ such that
$$
\left\|\Psi_t(f)(\lambda) -  \sum_{i,j \in [-K,K]} \varphi_t\left( \tau_{t} \left(\alpha_i^{\lambda} \alpha_j^{\lambda}\right)f\right) \otimes e_{ij} \right\| \leq \epsilon
$$
for all $\lambda$. Since $\lambda \mapsto  \sum_{i,j \in [-K,K]} \varphi_t\left( \tau_{t} \left(\alpha_i^{\lambda} \alpha_j^{\lambda}\right)f\right) \otimes e_{ij}$ clearly is an element of $IB\otimes \mathbb K$, we conclude that $\Psi_t(f) \in IB \otimes \mathbb K$ for all $t \in \mathbb R$. Since the above approximation to $\Psi_t(f)$ can be made uniform in any compact subset of $\mathbb R$, we see that $t \mapsto \Psi_t(f)$ is continuous. We claim that $\Psi = \left(\Psi_t\right)_{t \in \mathbb R}: SA \to IB$ is asymptotic homomorphism. It follows from Theorem \ref{aA} that 
$$
\left\{ \tau_t\left(\alpha_i^{\lambda}\alpha_j^{\lambda} \right) f: \ t \in \mathbb R, i,j \in \mathbb Z, \lambda \in [0,1] \right\}
$$
is a pre-compact subset of $SA$. We can therefore conclude from i) Lemma \ref{pre} that $\lim_{t \to -\infty} \Psi_t(f) = 0$ for all $f \in SA$. If we apply iii) and iv) of Lemma \ref{aA} instead of i), and use (\ref{uu13}) we find that $\lim_{t \to \infty} \Psi_t(f)\Psi_t(g) - \Psi_t(fg) = 0$ for all $f,g \in SA$, and similar considerations regarding linearity and self-adjointness show that $\Psi$ is an asymptotic homomorphism. Similarly, an obvious application of v) of Lemma \ref{pre} shows that $\Psi$ is asymptotically translation invariant. Thus $\Psi$ gives us a homotopy of asymptotically translation invariant asymptotic homomorphisms connecting $\psi''$ to $s \circ \varphi$ and we have established the commutativity of the first triangle. To obtain it for the second, it suffices to note that $\Psi$ is translation invariant when $\varphi$ is; the algebraic condition, $\Psi_{t-s} = \Psi_t \circ \tau_s$, is trivial, and the equi-continuity follows from the same considerations we presented for $\Phi$ above.
\end{proof}

By standard arguments $s_* : \Ext_h(A,B) \to \Ext_h(A,B \otimes \mathbb K)$ and $s_* : [[SA,B]]_{a, \tau} \to [[SA, B \otimes \mathbb K]]_{a, \tau}$ are both isomorphisms (of semigroups) when $B$ is stable, so Lemma \ref{l10} and Lemma \ref{l12} yields the following:

\begin{thm}\label{th3}  Let $A$ and $B$ be $C^*$-algebras, $A$ separable, $B$ $\sigma$-unital and stable. Then $[[SA,B]]_{a, \tau} \simeq  [[SA,B]]_{\tau} $ as abelian semigroups, and
$$
CH_{\tau} : \Ext_h(A,B) \to [[SA,B]]_{ \tau} 
$$
is an isomorphism.
\end{thm}
\begin{proof} It follows from Lemma \ref{l10} and Lemma \ref{l12} that $CH_{\tau} : \Ext_h(A,B) \to [[SA,B]]_{a,\tau}$ is an isomorphism. It follows from Lemma \ref{l12} that the diagram 
\begin{equation*}
\begin{xymatrix}{
&\Ext_h(A,B\otimes\mathbb K)\ar[dr]_-{CH_\tau}&\\
[[SA,B]]_{a,\tau}\ar[ur]_-{I} \ar[drr]_(0.67){s_*}|\hole&&
[[SA,B\otimes\mathbb K]]_\tau\ar[d]\\
[[SA,B]]_\tau\ar[u] \ar[urr]_(0.3){s_*}&&[[SA,B\otimes\mathbb K]]_{a,\tau} } 
\end{xymatrix}
\end{equation*}
commutes. Since the two $s_*$-maps are both isomorphisms, we conclude that the map $[[SA,B]]_{\tau} \to [[SA,B]]_{a,\tau}$ also is. 
\end{proof}

It should be pointed out that translation invariance occurs in other forms in other (related) constructions: In \cite{H} Higson bases a proof of the Atiyah-Singer index theorem on the construction of an asymptotic homomorphism $T = \left(T^{\omega}\right)_{\omega \in [1,\infty)} : C_0(T^*M) \to \mathbb K\left(L^2(M)\right)$, where $M$ is a smooth manifold without boundary, equipped with a smooth measure, and $T^*M$ is the cotangent bundle of $M$. Since every vector space carries a canonical action $\tau$ by the multiplicative group $\mathbb R^+$, viz. $\tau_s(x) = sx$, so does $T^*M$ and in turn also $C_0(T^*M)$. It is easy to see that Higson's asymptotic homomorphism is asymptotically translation invariant in the sense that
$$
\lim_{\omega \to \infty} \left\|T^{\omega s}(a) - T^{\omega}\tau_s(a)\right\| = 0
$$
for all $a \in C_0(T^*M)$ and all $\mathbb R^+$. The same is the case of many of the asymptotic homomorphisms constructed in \cite{HKT}.

As a final remark, let us point out that $[[SA,B \otimes \mathbb K]]_{\tau}$ is known to coincide with the E-theory group $E(SA,B)$ in at least two cases: When $A$ is nuclear and when $A$ is a suspension. It remains unclear exactly how widespread this coincidence is.

\end{document}